\documentclass{article}
\usepackage{authblk}
\usepackage{amscd}
\usepackage{graphicx}
\usepackage{amsmath}
\usepackage{amssymb}
\usepackage{amsthm}
\usepackage{cite}

\newtheorem{dfn}{Definition}[section]
\newtheorem{thm}[dfn]{Theorem}
\newtheorem{prop}[dfn]{Proposition}
\newtheorem{lem}[dfn]{Lemma}

\newtheorem{rem}[dfn]{Remark}

\newtheorem{prob}[dfn]{Problem}

\numberwithin{equation}{section}

\usepackage{color}

\begin{document}

\title{Characterization of Maximizers in A Non-Convex Geometric Optimization Problem With Application to Optical Wireless Power Transfer Systems}

\author{Dinh Hoa Nguyen \thanks{International Institute for Carbon-Neutral Energy Research (WPI-I2CNER), and Institute of Mathematics for Industry (IMI), Kyushu University, Fukuoka 819-0395, Japan, {\tt hoa.nd@i2cner.kyushu-u.ac.jp}. Corresponding author.}, ~
Kaname\, Matsue \thanks{Institute of Mathematics for Industry / International Institute for Carbon-Neutral Energy Research (WPI-I2CNER), Kyushu University, Fukuoka 819-0395, Japan, {\tt kmatsue@imi.kyushu-u.ac.jp}}}


\date{}
\maketitle


\begin{abstract}
This research studies a non-convex geometric optimization problem arising from the field of optical wireless power transfer. In the considered optimization problem, the cost function is a sum of negatively and fractionally powered distances from given points arbitrarily located in a plane to another point belonging to a different plane. Therefore, it is a strongly nonlinear and non-convex programming, hence posing a challenge on the characterization of its optimizer set, especially its set of global optimizers. To tackle this challenge, the bifurcation theory is employed to investigate the continuation and bifurcation structures of the Hessian matrix of the cost function. As such, two main results are derived. First, there is a critical distance between the two considered planes such that beyond which a unique global optimizer exists. Second, the exact number of maximizers is locally derived by the number of bifurcation branches determined via one-dimensional isotropic subgroups of a Lie group acting on $\mathbb{R}^2$, when the inter-plane distance is smaller than the above-mentioned critical distance. Consequently, numerical simulations and computations of bifurcation points are carried out for various configurations of the given points, whose results confirm the derived theoretical outcomes. 
\end{abstract}

\section{Introduction}

Optical wireless power transfer (OWPT) has recently gained much attention as an emerging and promising technology, due to its appealing characteristics over other wireless power transfer (WPT) technologies (e.g., inductive WPT or capacitive WPT), for instance longer transmitting distances and robustness to electromagnetic interferences \cite{9318756,Nguyen-fenr-21,9422160}. OWPT can be employed for systems and applications in internet of things \cite{8398229,8618313}, consumer electronics \cite{8398229,7456176}, implantable biomedical devices \cite{Kim2020ActivePW,JEONG2021113139,Lu21}, vehicle electrification \cite{Nguyen-fenr-21,Nguyen-u2wpt,Nguyen-Access20,Nguyen-Access23,Manson11}, robotics \cite{8970169,KJin19}, and so on. Moreover, OWPT can be used in different environments such as air, water and undersea \cite{8579163,8970169}, under clothing \cite{su15108146}, etc. 

For OWPT systems, there are two typical types of optical transmitters, namely laser and light emitting diode (LED). Laser transmitters have high intensity and well-focused output lights, hence being able to wirelessly deliver more power and over a longer distance than LEDs \cite{8404104,1360294812944463744,9903468,MOHAMMADNIA2021107283,Gou:22}. Nevertheless, laser-based OWPT incurs higher complexity, larger size, and higher cost for the system, due to many system components need to be added for generating laser beams and ensuring their safety \cite{9422160}. On the other hand, LED transmitters have a stronger tolerance to misalignments with the optical receivers, due to their spread output lights. Moreover, LEDs are safer to users, because of their low intensity outputs. Therefore, LED-based OWPT systems have been investigated in several recent works, e.g. \cite{photonics10070824,photonics9080576,Nguyen-photonics-22,Nguyen-fenr-21,9081113}. 

It is worth emphasizing that most of existing works on LED-based OWPT considered a single LED as the optical transmitter, while other few studies on OWPT systems using LED arrays, see e.g. \cite{photonics10070824,photonics9080576}, did not derive a mathematical expression to represent the power transmitting efficiency from the transmitter to the receiver. A recent work \cite{su15108146} has proposed a mathematical formula aiming to close the above-mentioned research gap, in which the output power of the solar cell optical receiver is linearly approximated by the summation of that obtained from individual LEDs in an LED array. 

Interestingly, the mathematical formula in \cite{su15108146} exhibits rich behaviors as the transmitting distance is changed, depicted through simulation results. For example, if the transmitting distance is longer than a certain value, then there exists a single position for the optical receiver at which the received wireless power is maximum. However, at transmitting distances smaller than such value, there could be more than one position at which the received wireless power is maximum. 

In the language of mathematical programming, the descriptions above mean that there could be different numbers of maximizers for the optimization problem whose cost function is to maximize the received wireless power transmitted from a bunch of LEDs, depending on the transmitter-receiver distance. These complicated features of the maximizer set is due to the non-convexity nature of the considered optimization problem, as will be clearly seen in the problem description in Section \ref{prob}.  

Motivated by the realistic system of OWPT using multiple LEDs shown above, this research attempts to investigate the characteristics of the optimizer set of a  geometric optimization problem which is more general but close to that for the multi-LED based OWPT systems. More specifically, the geometric optimization problem in the current work considers the maximization of a sum of negatively and fractionally powered distances from given points arbitrarily located in a plane to another point belonging to a different plane (see Fig. \ref{OWPT_example} for an illustration). Then our aim is to figure out conditions under which a single maximizer exists and the switch from such a single global maximizer to multiple global/local maximizers occurs, as well as mathematical foundations underlying the existence of such multiple global/local maximizers.

Non-convex and geometric optimization problems have been actively studied in the literature. The work in \cite{zhang2022symmetry} presented a survey on a class of tractable non-convex optimization problems with rotational and discrete symmetries which can be found in machine learning and data-driven applications. Likewise, \cite{MAL-058} introduced a review of non-convex optimization problems and several solving techniques, with the focus on the machine learning applications. It can be seen from those reviews that a rich body of works has been made on low-rank, sparse matrix factorization or robust linear regression problems. Nevertheless, geometric optimization problems have been less investigated. 
In \cite{140978168}, geometric optimization problems with geodesically convex and log-nonexpansive cost functions were considered. The non-uniqueness and symmetry of solution in a topology optimization problem were studied in \cite{Watada10}. On the other hand, the research in \cite{Levin24} presented a framework to study the parametrization effect on the local minima and critical points of non-convex cost functions.  
A tutorial on the non-convex optimization problems related to geometry and manifolds, and the methods to solve them can be found in \cite{9194028}. An interesting example given in \cite{9194028} was the Riemannian center of mass which is defined as the minimizer of a function equal to a sum of squared geodesic distances between a point to a set of given points on a Riemann manifold. This problem is close to our considering geometric non-convex optimization problem described above, but is much more simple where it is reduced to the conventional average distance on a Euclidean space, whereas our considering problem has no simple and analytical solution like that.  

It is also noteworthy that very few results have been existed so far on the characterization of the number of global/local optimizers for non-convex, nonlinear optimization problems including geometric mathematical programming. Therefore, the current paper aims to fulfill this research gap by explicitly revealing a mathematical framework to determine exactly the number of maximizers in the considered non-convex, nonlinear geometric optimization problem aforementioned. 
Our specific contributions are as follows.
\begin{itemize}
	\item We prove the existence of a positive {\it lower bound} for the distance between a moving point and the plane hosting given points such that the considering optimization problem has a {\it unique global maximizer} whenever the distance is greater than this bound. This is held for arbitrary locations of given points. 
	\item We establish {\it a connection between the fields of non-convex geometric optimization and bifurcation theory}, where the latter lays the foundation for characterizing the exact number of local maximizers in the considering non-convex geometric optimization problem. More specifically, {\it the exact number of maximizers is locally derived by the number of bifurcation branches determined via one-dimensional isotropic subgroups of a Lie group acting on $\mathbb{R}^2$}. To the best of our knowledge, this is the first time such result is reported.   
\end{itemize}

The remainder of this paper is as follows. The description of the considered geometric optimization problem and details on the motivating example of OWPT are given in Section \ref{prob}. Then the main theoretical results are presented in Section \ref{results}. Simulation and numerical computation results to illustrate the derived theoretical results are provided in Section \ref{num}. Finally, conclusions and directions for future works are included in Section \ref{concl}.

\section{Problem Setting}
\label{prob}

\subsection{Problem Description}

Let $\Pi_1, \Pi_2$ be $2$-dimensional subspaces in $\mathbb{R}^3$ subspaces defined by
\begin{equation*}
\Pi_1 := \{z = 0\}, \quad \Pi_2 := \{z = h\}
\end{equation*}
for given $h > 0$.
Let ${\bf p}_1, \ldots, {\bf p}_n$ be given points in $\Pi_1$ with the coordinate ${\bf p}_i = ({\bf x}_i, 0) \equiv (x_i, y_i, 0)$.
The main problem we shall consider is the following.

\begin{prob}
\label{prob-def}
Let ${\bf p} = ({\bf x}, h)\equiv (x,y, h) \in \Pi_2$. 
Then, for given points $\{{\bf p}_i\}_{i=1}^n \subset \Pi_1$ and $h > 0$, determine the point ${\bf p}$ maximizing the following functional:
\begin{align*}
&f({\bf p}) \equiv f({\bf x}, h) = f({\bf x}, h; {\bf x}_1, \ldots, {\bf x}_n) := \sum_{i=1}^n d_i^{-m},\\
&d_i := \left\{ (x-x_i)^2 + (y - y_i)^2 + h^2 \right\}^{1/2},\quad m\in \left(3, 4\right).
\end{align*}
\end{prob}

In the following expression, the point ${\bf p}$ is always assumed to be in $\Pi_2$ and hence we identify the expression of ${\bf p} = (x,y,h)$ with ${\bf x} = (x,y)$.
Similarly, given points ${\bf p}_i = (x_i, y_i, 0)$, $i=1,\ldots, n$, are identified with ${\bf x}_i = (x_i, y_i)$.


In the present study, we interpret the maximizer ${\bf p} = ({\bf x}, h)$ as the critical point of the functional $f$ with fixed $h$, which solves
\begin{equation*}
F({\bf x}, h) = (F_1({\bf x}, h), F_2({\bf x}, h) )\equiv \nabla_{\bf x} f({\bf x}, h) := \left( \frac{\partial f}{\partial x}({\bf x}, h), \frac{\partial f}{\partial y}({\bf x}, h)\right)^T = 0.
\end{equation*}
When we are interested in the local (or global) maximizer, it is sufficient to verify if the Jacobian matrix 
\begin{equation*}
D_{\bf x}F({\bf x}, h) \equiv \begin{pmatrix}
f_{xx}({\bf x}, h) & f_{xy}({\bf x}, h)\\
f_{yx}({\bf x}, h) & f_{yy}({\bf x}, h)
\end{pmatrix},
\end{equation*}
namely the Hessian matrix of $f$ at the maximizer ${\bf p} = ({\bf x}, h)$, is negative-definite.


For any {\em fixed} $\{{\bf x}_i\}_{i=1}^n$ and $h > 0$, the problem is a $2$-dimensional zero-finding problem with $(2n+1)$-parameters.
Hence the continuation and bifurcation structure of zeros can be computed very rapidly for specified parameters.
In particular, numerical bifurcation theory can be applied very effectively to computing parameter families of (local) maximizers.

\subsection{Motivating Example}

\begin{figure}[htpb!]
		\centering
		\includegraphics[scale=0.5]{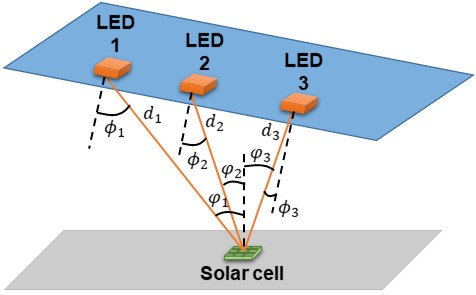}
		\caption{Illustration for an optical wireless power system consisting of a multi-LED transmitter.}
		\label{OWPT_example}
	\end{figure}

Consider a practical system of optical wireless power and data transfer (OWPDT) using light emitting diodes (LED) and solar cells, as illustrated in Fig. \ref{OWPT_example}. In this system, light from an LED array, for example that located on the room ceiling, is modulated to illuminate on a device equipped with a solar cell so that power, data, or both can be wirelessly sent from such an LED array to the device. This provides an innovative wireless powering and information exchange approach using visible light that possesses several advantages over other wireless powering and information exchange methods. First, it is known that such OWPDT system based on visible light is more robust to interferences than those based on radio or microwave signals, hence it has been considered for the next-generation, i.e. 6G communication networks. Second, this type of OWPDT system is safer for human and other living objects than those utilizing laser or other electromagnetic waves. Last but not least, this LED-based OWPDT system is reasonably simple and low-cost where the existing LED lighting system could be inherited and supplemented. Therefore, LED-based OWPDT systems has been an emerging research area recently.   

For the above-illustrated LED-based OWPDT system, assume that line-of-sight (LOS) optical links can be established between $n$ LEDs and the solar cell, then the LOS DC gain of each optical link between each LED and the solar cell can be estimated by,
\begin{equation}
\label{los-dc-gain}
\eta = \left\{ 
\begin{aligned}
	\dfrac{(m_l+1)A}{2\pi d_i^2} \cos^{m_l}\phi_i \cos\varphi_i: & ~0 \leq \varphi_i \leq \bar{\varphi}, \\
	0: & ~\textrm{otherwise}.
\end{aligned}
\right.
\end{equation}
The variables in Eq. \eqref{los-dc-gain} are as follows. $A$ is the physical area of the solar cell, $d_i$ is the distance between the $i$th LED to the solar cell, $\phi_i$ is the angle of incidence (i.e. the angle at which the receiver sees the transmitter), $\varphi_i$ is the angle of irradiance (i.e. the angle at which the transmitter sees the receiver), $\bar{\varphi}$ is the width of the field of view (FOV) at the optical receiver. In addition, $m_l$ denotes the order of Lambertian emission which is given by the semi-angle at half illuminance $\Phi_{1/2}$ of an LED, where $m_l = -\dfrac{\ln 2}{\ln \cos\Phi_{1/2}}$. 

Note that $\cos\phi_i$ and $\cos\varphi_i$ are positive and smaller as $\phi_i$ and $\varphi_i$ are larger, as long as they are in the interval $(0,\pi/2)$. Hence, $\eta$ can be further simplified by only the fist equation in \eqref{los-dc-gain}, without the constraint. As such, the overall LOS DC gain made by the LED array can be approximated by,
\begin{equation}
\label{gain-total}
\eta_{\textrm{total}} = \dfrac{(m_l+1)A}{2\pi} \sum_{i=1}^{n} \dfrac{1}{d_i^2} \cos^{m_l}\phi_i \cos\varphi_i,
\end{equation}
assuming that the order of Lambertain $m_l$ of all LEDs are the same. 

Consequently, suppose that the planes containing the LEDs and the solar cell are parallel, and the distance between is denoted by $h$. Then $\phi_i = \varphi_i$ for all $i=1,\ldots,n$. Furthermore, $\cos\phi_i = \dfrac{h}{d_i}, i=1,\ldots,n$. Thus, Eq.\eqref{gain-total} becomes,
\begin{equation}
\label{gain-total-1}
\eta_{\textrm{total}} = \dfrac{(m_l+1)A h^{m_l+1}}{2\pi} \sum_{i=1}^{n} \dfrac{1}{d_i^{m_l+3}}. 
\end{equation} 
Omitting the constant $\dfrac{(m_l+1)A h^{m_l+1}}{2\pi}$, the rest in Eq. \eqref{gain-total-1} is exactly in the form of $f({\bf p})$ with $m=m_l+3$, since $0 < m_l <1$. 
This practical and meaningful system of OWPDT motivates us to consider the general problem specified in the problem \eqref{prob-def}.

\section{Main Results}
\label{results}

\subsection{Characterization for the scenario of a unique global maximizer}

Our first main theoretical result is given in the following theorem.

\begin{thm}
\label{thm-existence-max}
For {\bf any} configurations of given points $\mathcal{X}_n \equiv \{{\bf x}_i\}_{i=1}^n$, there is a positive number $h_0 = h_0(\mathcal{X}_n) > 0$ such that, for {\bf any} $h > h_0$, the functional $f({\bf x}, h)$ admits {\bf a unique global maximizer} ${\bf p} = {\bf p}(h, \mathcal{X}_n)$.
\end{thm}

First observe that, for any $h > 0$ as a parameter, $f({\bf x}, h)$ as a function of ${\bf x} = (x,y)$ is smooth on $\mathbb{R}^2$.
We therefore know that, for any local maximizer or minimizer ${\bf x}_\ast = (x_\ast, y_\ast)$ of $f$ must be critical points, namely
\begin{equation*}
\nabla_{\bf x} f({\bf x}, h) \equiv F({\bf x}, h) = \begin{pmatrix}
F_1({\bf x}, h) \\ F_2({\bf x}, h)
\end{pmatrix} = \begin{pmatrix}
0 \\ 0
\end{pmatrix}\quad \text{ at }\quad {\bf x} = {\bf x}_\ast. 
\end{equation*}
Our strategy for the proof of Theorem \ref{thm-existence-max} is decomposed into the following steps.
\begin{enumerate}
\item Maximizers or minimizers of $f$ must be located in a compact subset $R$ of $\mathbb{R}^2$.
As a result, $f$ admits at least one maximizer and minimizer.
\item The Jacobian matrix $D_{\bf x}F({\bf x},h)$ of $F$, namely the Hessian matrix 
of $f$ with respect to ${\bf x}$, is negative definite for all ${\bf x}\in R$, provided $h$ is sufficiently large. 
In particular, $f$ is strictly concave on $R$.
\end{enumerate}
Before starting the proof, observe that the gradient $F = \nabla_{\bf x} f$ at ${\bf p}=({\bf x},h)$ is given by
\begin{equation}
\label{grad-f}
F({\bf x}, h) = \begin{pmatrix}
F_1({\bf x}, h) \\ F_2({\bf x}, h)
\end{pmatrix} = \left( -m \sum_{i=1}^n \frac{x-x_i}{d_i^{m+2}},\quad -m \sum_{i=1}^n \frac{y-y_i}{d_i^{m+2}}\right)^T.
\end{equation}

\subsubsection{Step 1 in the proof of Theorem \ref{thm-existence-max}}
 
\begin{lem}
\label{lem-max-location}
Maximizers and minimizers of $f$ are located in the rectangular domain 
\begin{equation*}
R = R(\mathcal{X}_N) := \left[ \min_{i=1,\ldots, n} x_i, \max_{i=1,\ldots, n} x_i\right] \times \left[ \min_{i=1,\ldots, n} y_i, \max_{i=1,\ldots, n} y_i \right].
\end{equation*}
Moreover, $f$ does not attain a maximal value outside $R$.
\end{lem}
\begin{proof}
Let ${\bf x} = (x,y) \in \mathbb{R}^2\setminus R$.
Then, at least one the following cases must be held: (i) $x < x_i$ for all $i$, (ii) $x > x_i$ for all $i$, (iii) $y < y_i$ for all $i$, or (iv) $y > y_i$ for all $i$.
\par
When (i) holds, it follows from (\ref{grad-f}) that $F_1({\bf x},h) > 0$. 
Similarly, $F_1({\bf x},h) < 0$ must be satisfied when (ii) holds, while $F_2({\bf x},h) > 0$ when (iii) holds, and $F_2({\bf x},h) < 0$ when (iv) holds, respectively.
Because all local maximizers and minimizers must be critical points of $f$, namely zeros of $F$, the above assertions imply that any point ${\bf x} \in \mathbb{R}^2\setminus R$ must be neither a maximizer nor a minimizer of $f$.
\par
Finally, observations in cases (i) and (ii) also imply that $f$ must increase as $x$ moves toward
\begin{equation*}
R_x \equiv \left[ \min_{i=1,\ldots, n} x_i, \max_{i=1,\ldots, n} x_i\right],
\end{equation*}
the $x$-component of $R$.
The similar assertion holds for $y$, indicating that the maximizer ${\bf x}$ of $f$, if it exists, must be in $R$ and the proof is completed.
\end{proof}

Because our aim here is to find the global maximizer of $f$, we pay our attention to maximizers, which have to be in $R$.
Thanks to the smoothness of $f$ and the compactness of $R$, it admits at least one (local) maximizer in $R$.
The rest of our discussions here is the uniqueness.

\subsubsection{Step 2 in the proof of Theorem \ref{thm-existence-max}}

\begin{prop}
\label{prop-neg-def}
For any compact sets $K\subset \mathbb{R}^2$, there is a positive number $h_0 =h_0(K) > 0$ such that, for all $h > h_0$, $D_{\bf x}F({\bf x},h)$ is negative definite for all ${\bf x}\in K$.
\end{prop}

\begin{proof}
It is sufficient to prove that both eigenvalues of $D_{\bf x}F({\bf x},h)$ are negative for all ${\bf x} = (x,y)$ with a suitable choice of $h$.
Now the Jacobian matrix $D_{\bf x}F({\bf x},h)$ is given by $D_{\bf x}F({\bf x},h) = \begin{pmatrix} a & b \\ b & c \end{pmatrix}$, where
\begin{align*}
a &= -m  \sum_{i=1}^n \frac{1}{d_i^{m+2}}\left[ 1 - (m+2) \frac{(x - x_i)^2}{d_i^2}\right],\\ 
b &= m(m+2) \sum_{i=1}^n \frac{(x - x_i)(y - y_i)}{d_i^{m+4}}, \quad
c = -m  \sum_{i=1}^n \frac{1}{d_i^{m+2}}\left[ 1 - (m+2) \frac{(y - y_i)^2}{d_i^2}\right].
\end{align*}
Because $D_{\bf x}F(x,y)$ is symmetric, both eigenvalues are real and hence it is sufficient to prove $a + c< 0$ and $ac > b^2$ for our statement.
\par
First observe that
\begin{align*}
a +c &= -m  \sum_{i=1}^n \frac{1}{d_i^{m+2}}\left[ 1 - (m+2) \frac{(x - x_i)^2}{d_i^2}\right] - m  \sum_{i=1}^n \frac{1}{d_i^{m+2}}\left[ 1 - (m+2) \frac{(y - y_i)^2}{d_i^2}\right] \\
	&= -m\sum_{i=1}^n \frac{1}{d_i^{m+2}} \left( 2 - (m+2) \frac{(x - x_i)^2 + (y - y_i)^2}{d_i^2}  \right).
\end{align*}
We easily know that the rightmost function is negative whenever
\begin{equation}
\label{ineq-1-i}
\frac{(x - x_i)^2 + (y - y_i)^2}{(x - x_i)^2 + (y - y_i)^2 + h^2} < \frac{2}{m+2}
\end{equation}
holds for all $i$.
For fixed ${\bf x}\in K$, we can choose $\tilde h_{0,i}^1 = \tilde h_{0,i}^1({\bf x})$ such that (\ref{ineq-1-i}) for all $h > \tilde h_{0,i}^1$.
Choosing $\tilde h_0^1({\bf x}) := \max_{i=1,\ldots, n}\tilde h_{0,i}^1({\bf x})$, the inequality (\ref{ineq-1-i}) holds for {\em all $i$} and all $h > \tilde h_0^1$, with fixed ${\bf x}\in K$.
By continuity of $a+c$ as a function of ${\bf x}$, we can choose an open neighborhood $U_{{\bf x}}\subset \mathbb{R}^2$ of ${\bf x}$ such that (\ref{ineq-1-i}) holds for all $i$, all $h > \tilde h_0^1$ and all $\tilde {\bf x} = (\tilde x,\tilde y)\in U_{{\bf x}}$.
By compactness of $K$, there is a finite subcovering $\{U_{{\bf x}_j'}\}_{j=1}^l$ of $K$ with $\{{\bf x}_j'\}_{j=1}^l \subset K$.
Therefore, letting
\begin{equation*}
h_0^1 := \max_{j=1,\ldots, l} \tilde h_0^1({\bf x}_j'),
\end{equation*}
it follows that the inequality (\ref{ineq-1-i}) holds for all $i$, all $h > h_0^1$ and {\em all ${\bf x}\in K$}.
\par
Next, observe that
\begin{align}
	\label{neg-2}
& ac - b^2 \notag \\ 
 &= \left(   \sum_{i=1}^n \frac{-m}{d_i^{m+2}}\left[ 1 - (m+2) \frac{(x - x_i)^2}{d_i^2}\right]\right) \left(  \sum_{i=1}^n \frac{ -m}{d_i^{m+2}}\left[ 1 - (m+2) \frac{(y - y_i)^2}{d_i^2}\right] \right) \notag \\
&\qquad - m^2(m+2)^2 \left(\sum_{i=1}^n \frac{(x - x_i)(y - y_i)}{d_i^{m+4}}\right)^2 \notag \\
	&= m^2\left\{ \left( \sum_{i=1}^n \frac{1}{d_i^{m+2}}\left[ 1 - (m+2) \frac{(x - x_i)^2}{d_i^2}\right]\right) \left( \sum_{i=1}^n \frac{1}{d_i^{m+2}}\left[ 1 - (m+2) \frac{(y - y_i)^2}{d_i^2}\right] \right) \right. \notag \\ 	
&\qquad - \left. (m+2)^2 \left(\sum_{i=1}^n \frac{(x - x_i)(y - y_i)}{d_i^{m+4}}\right)^2 \right\}.
\end{align}
We shall estimate the rightmost function in (\ref{neg-2}) component-wise.
First we have
\begin{align}
\label{neg-2-1}
&\frac{1}{d_i^{2(m+2)}}\left[ 1 - (m+2) \frac{(x - x_i)^2}{d_i^2}\right]\cdot \left[ 1 - (m+2) \frac{(y - y_i)^2}{d_i^2}\right] \notag \\
&- (m+2)^2\frac{(x - x_i)^2(y - y_i)^2}{d_i^{2(m+4)}} \nonumber \\
&= \frac{1}{d_i^{2(m+2)}} \left\{ 1 - (m+2) \frac{(x - x_i)^2 + (y - y_i)^2}{d_i^2}\right\},
\end{align}
and the right-hand side is estimated later.
Next, consider
\begin{align*}
& \frac{1}{d_i^{m+2}d_j^{m+2}}\left[ \left\{ 1 - (m+2) \frac{(x - x_i)^2}{d_i^2}\right\} \left\{ 1 - (m+2) \frac{(y - y_j)^2}{d_j^2}\right\} \right. \\
& + \left. \left\{ 1 - (m+2) \frac{(x - x_j)^2}{d_j^2}\right\} \left\{ 1 - (m+2) \frac{(y - y_i)^2}{d_i^2} \right\} \right] \\
& - 2(m+2)^2\frac{(x - x_i)(y - y_i)}{d_i^{m+4}}\frac{(x - x_j)(y - y_j)}{d_j^{m+4}}
\end{align*}
with $i\not = j$.
Observe that
\begin{align*}
&\left\{ \frac{(x - x_i)^2}{d_i^2} \frac{(y - y_j)^2}{d_j^2} + \frac{(x - x_j)^2}{d_j^2} \frac{(y - y_i)^2}{d_i^2} \right\} - 2 \frac{(x - x_i)(y - y_i)}{d_i^2}\frac{(x - x_j)(y - y_j)}{d_j^2} \\
 &= \frac{1}{d_i^2d_j^2}\left\{ (x - x_i)^2 (y - y_j)^2 + (x - x_j)^2(y - y_i)^2 - 2(x - x_i)(y - y_i)(x - x_j)(y - y_j) \right\} \\
 &= \frac{1}{d_i^2d_j^2}\left\{ (x - x_i) (y - y_j) - (x - x_j)(y - y_i) \right\}^2 \geq 0.
\end{align*}
Then the rest is to estimate
\begin{align}
\label{neg-2-2}
2 - (m+2) \frac{(x - x_i)^2 + (y - y_i)^2}{d_i^2} - (m+2) \frac{(x - x_j)^2 + (y - y_j)^2}{d_j^2}.
\end{align}
Here it is concluded that both (\ref{neg-2-1}) and (\ref{neg-2-2}) are always satisfied whenever
\begin{equation}
\label{neg-2-3}
\frac{(x - x_i)^2 + (y - y_i)^2}{(x - x_i)^2 + (y - y_i)^2 + h^2} < \frac{1}{m+2},
\end{equation}
holds for all $i$.
The similar argument involving (\ref{ineq-1-i}) yields that the inequality (\ref{neg-2}) is true for all $h > \tilde h_{0,i}^2({\bf x})$ with some $\tilde h_{0,i}^2({\bf x}) > 0$.
Also, by the similar arguments to obtaining $h_0^1$, we know that there is a number $h_0^2 > 0$ such that (\ref{neg-2}) holds for all $i$, all $h > h_0^1$ and all ${\bf x}\in K$.
\par
Finally, choosing $h_0 := \max\{h_0^1, h_0^2\}$, we conclude that both inequalities $a + c < 0$ and $ac > b^2$ holds for all ${\bf x} \in K$ and all $h > h_0$, which proves our claim.
\end{proof}

Choosing $h_0$ with $K := R$ as in Proposition \ref{prop-neg-def} with $R$ defined in Lemma \ref{lem-max-location}, we know that:
\begin{itemize}
\item $f$ has a (local) maximizer in $R$, and is strictly concave in $R$ in the sense that the Hessian matrix ${\rm Hess}(f)({\bf x})$ is negative definite for all ${\bf x}\in R$, whenever $h > h_0$,
\item $f$ does not attain maximal values outside $R$.
\end{itemize}
As a conclusion, the functional $f$ admits a unique global maximizer in $R$ and the proof of Theorem \ref{thm-existence-max} is completed.

\subsection{Symmetry of functionals and possible maximizers}
\label{section-symmetry}

The previous arguments guarantee the {\em unique} existence of maximizers of $f$ when $h > 0$ is relatively large.
As seen below, there are possibilities to exist multiple (local) maximizers of $f$ when $h$ becomes considerably small.
In typical situations, the existence of multiple maximizers is a consequence of {\em bifurcations} of solutions to $F({\bf x},h) = 0$ and those bifurcations can be extracted through standard machineries in (numerical) bifurcation theory (e.g., \cite{K2013}).
On the other hand, when there are {\em symmetries} in our objects, multiple maximizers and/or their bifurcations can be considered. 
The simplest example to easily observe symmetry in our considering problem is the configuration of given points, such as squares, hexahedra or circular distributions.
They possess reflection symmetry $(x, y)\mapsto (\pm x, \pm y)$ and/or rotations with fixed angle $\theta$:
\begin{equation}
\label{rotation}
\begin{pmatrix}
x \\ y
\end{pmatrix}\mapsto \begin{pmatrix}
\cos \theta & -\sin \theta \\ 
\sin \theta & \cos \theta
\end{pmatrix}\begin{pmatrix}
x \\ y
\end{pmatrix}.
\end{equation}
In addition to the configuration of given points, the functional $F = (F_1, F_2)$ should be paid attention to for extracting symmetry of maximizers.
Symmetry of $F$ can be restricted compared with that of the given points' configuration. 
Here we discuss the possible symmetry on $F$ following the theory of {\em bifurcations with symmetry} (see e.g., \cite{GSS1988}).

Mathematically, symmetries on variables and functionals are described through {\em groups}. 
Let $\Gamma$ be a (compact Lie) group\footnote{
Our examples of $\Gamma$ are $\mathbb{Z}_n$ (cyclic group) or $D_n$ (dihedral groups) mentioned below.
Readers should not be worried much about the details of Lie groups.
}.
The {\em (linear) action of $\Gamma$ on the vector space $V$} is the mapping 
\begin{equation*}
\Gamma\times V\to V,\quad (\gamma, {\bf v}) \mapsto \gamma\cdot {\bf v}\quad ({\bf v}\in V)
\end{equation*}
such that:
\begin{enumerate}
\item For each $\gamma \in \Gamma$ the mapping $\rho_\gamma: V\to V$ defined by $\rho_\gamma({\bf v}) = \gamma\cdot {\bf v}$ is linear.
\item If $\gamma_1, \gamma_2 \in \Gamma$, then $\gamma_1\cdot (\gamma_2\cdot {\bf v}) = (\gamma_1\gamma_2)\cdot {\bf v}$.
\end{enumerate}
Examples of compact Lie groups exhibiting such actions are {\em cyclic groups} $\mathbb{Z}_n$ or {\em dihedral groups} $D_n$ for some integer $n$.
Several actions of these groups we shall frequently use here are the following:
\begin{description}
\item[$\mathbb{Z}_2$-actions on $\mathbb{R}^2$.]
\end{description}
Either of them:
\begin{equation*}
\begin{pmatrix}
x \\ y
\end{pmatrix}\mapsto \begin{pmatrix}
-x \\ y
\end{pmatrix},\quad 
\begin{pmatrix}
x \\ y
\end{pmatrix}\mapsto \begin{pmatrix}
x \\ -y
\end{pmatrix},\quad
\begin{pmatrix}
x \\ y
\end{pmatrix}\mapsto \begin{pmatrix}
-x \\ -y
\end{pmatrix},
\end{equation*}
namely {\em reflection} in one axis or two axes simultaneously.
\begin{description}
\item[A $D_n$-action on $\mathbb{R}^2$.]
\end{description}
The following transformations and their compositions:
\begin{equation}
\label{Dn} 
\xi_n: \begin{pmatrix}
x \\ y
\end{pmatrix}\mapsto \begin{pmatrix}
\cos \frac{2\pi}{n} & -\sin \frac{2\pi}{n} \\ 
\sin \frac{2\pi}{n} & \cos \frac{2\pi}{n}
\end{pmatrix}\begin{pmatrix}
x \\ y
\end{pmatrix},\quad 
\kappa: \begin{pmatrix}
x \\ y
\end{pmatrix}\mapsto \begin{pmatrix}
x \\ -y
\end{pmatrix},
\end{equation}
namely reflection and {\em rotation}.

\begin{dfn}\rm
Let $\Gamma$ be a (compact Lie) group acting linearly on a vector space $V$, and $g: V\to V$ be a mapping.
We say that $g$ is {\em $\Gamma$-equivariant} if
\begin{equation*}
g(\gamma \cdot {\bf v}) = \gamma \cdot g({\bf v})
\end{equation*}
holds for all $\gamma \in \Gamma$ and ${\bf v}\in V$.
\end{dfn}

We concentrate on the most symmetric case in $\mathbb{R}^2$ with discretely distributed configurations; the {\em circular} configuration.
In the following proposition, labeling of the LED configuration is changed for easier expressions of points and their images of actions.
\begin{prop}
\label{prop-equivariant-LED}
Let $\mathcal{X}_n = \{{\bf x}_i\}_{i=0}^{n-1} \subset \mathbb{R}^2$ be the LED configuration evenly distributed on a circle with radius $r > 0$.
Define the $D_n$-action on $\mathbb{R}^2$ by (\ref{Dn}).
Then, for all $h>0$, the functional $F$ is $D_n$-equivariant in ${\bf x}$.
\end{prop}

\begin{proof}
Without the loss of generality, we may fix the position of each LED ${\bf x}_k$, $k=0,\ldots, n-1$ as
\begin{equation*}
{\bf x}_k \equiv (x_k, y_k)^T = \left( r\cos (kw_n), r\sin (kw_n) \right)^T,\quad w_n = \frac{2\pi}{n}.
\end{equation*}
Now the dihedral group $D_n$ generated by transformations (\ref{Dn}) is acting on $\mathbb{R}^2$, and the configuration $\mathcal{X}_n$ is $D_n$-invariant in the sense that, for any $\gamma\in D_n$ and ${\bf x}\in \mathcal{X}_n$, there is $\tilde {\bf x}\in \mathcal{X}_n$ such that $\gamma\cdot {\bf x} = \tilde {\bf x}$.
It is sufficient to show that 
\begin{equation*}
(F_1(\gamma\cdot {\bf x}, h), F_2(\gamma\cdot {\bf x}, h)) = \gamma\cdot (F_1({\bf x}, h), F_2({\bf x}, h))
\end{equation*}
for any $h > 0$ when $\gamma = \xi_n, \kappa$ given in (\ref{Dn}), thanks to the definition of the $D_n$-action on $\mathbb{R}^2$.
\par
First consider the action $\kappa$; the reflection across the $x$-axis.
Direct calculations yield that ${\bf x}\in \mathcal{X}_n$ if and only if $\kappa\cdot {\bf x}\in\mathcal{X}_n$.
More precisely, for ${\bf x} = {\bf x}_k$, the reflection $\kappa\cdot {\bf x}$ is given by ${\bf x}_{n-k}$ with the identification ${\bf x}_{n} = {\bf x}_0$.
We shall always use this identification unless otherwise noted.
Let ${\bf x} = (x,y)\in \mathbb{R}^2$ and consider the $\mathbb{Z}_2$-action given by
\begin{equation*}
\kappa: (x, y)^T \mapsto (x, -y)^T.
\end{equation*}
Then
\begin{align*}
F(x, -y, h) &= (F_1(x, -y, h), F_2(x, -y, h)),\\
\frac{-1}{m}F_1(x, -y, h) &= \sum_{i=0}^{n-1} \frac{x-x_i}{\{(x - x_i)^2 + (-y - y_i)^2 + h^2\}^{(m+2)/2}}\\
	&= \frac{x- r}{\{(x - r)^2 + y^2 + h^2\}^{(m+2)/2}} + \sum_{i=1}^{n-1} \frac{x-x_i}{\{(x - x_i)^2 + (-y - y_i)^2 + h^2\}^{(m+2)/2}}\\
	&= \frac{x- r}{\{(x - r)^2 + y^2 + h^2\}^{(m+2)/2}} + \sum_{i=1}^{n-1} \frac{x-x_i}{\{(x - x_i)^2 + (y + y_i)^2 + h^2\}^{(m+2)/2}}\\
	&= \frac{x- r}{\{(x - r)^2 + y^2 + h^2\}^{(m+2)/2}} + \sum_{i=1}^{n-1} \frac{x-x_{n-i}}{\{(x - x_{n-i})^2 + (y- y_{n-i})^2 + h^2\}^{(m+2)/2}}\\
	&= \sum_{i=0}^{n-1} \frac{x-x_i}{\{(x - x_i)^2 + (y - y_i)^2 + h^2\}^{(m+2)/2}},
\end{align*}
while
\begin{align*}
\frac{-1}{m}F_2(x, -y, h) &= \sum_{i=0}^{n-1} \frac{-y-y_i}{\{(x - x_i)^2 + (-y - y_i)^2 + h^2\}^{(m+2)/2}}\\
	&= \frac{-y}{\{(x - r)^2 + y^2 + h^2\}^{(m+2)/2}} + \sum_{i=1}^{n-1} \frac{-y-y_i}{\{(x - x_i)^2 + (-y - y_i)^2 + h^2\}^{(m+2)/2}}\\
	&= \frac{-y}{\{(x - r)^2 + y^2 + h^2\}^{(m+2)/2}} + \sum_{i=1}^{n-1} \frac{-y-y_i}{\{(x - x_i)^2 + (y + y_i)^2 + h^2\}^{(m+2)/2}}\\
	&= \frac{-y}{\{(x - r)^2 + y^2 + h^2\}^{(m+2)/2}} + \sum_{i=1}^{n-1} \frac{-y+y_{n-i}}{\{(x - x_{n-i})^2 + (y- y_{n-i})^2 + h^2\}^{(m+2)/2}}\\
	&= - \sum_{i=0}^{n-1} \frac{y-y_i}{\{(x - x_i)^2 + (y - y_i)^2 + h^2\}^{(m+2)/2}}.
\end{align*}
Consequently, we have $F(\kappa\cdot {\bf x}, h) = \kappa F({\bf x}, h) = (F_1({\bf x}, h), -F_2({\bf x}, h))$, i.e., $F$ is $\mathbb{Z}_2$-equivariant in ${\bf x}$.
Note that the group $\mathbb{Z}_2$ generated by $\kappa$ is a subgroup of $D_n$.
\par
Next, consider the action $\xi_n$; the $w_n$-rotation.
For given ${\bf x}= (x,y)^T\in \mathbb{R}^2$, let $\tilde {\bf x} = \xi_n\cdot {\bf x} \equiv (\tilde x, \tilde y)$.
More precisely,
\begin{equation*}
\tilde x = (\cos w_n)x - (\sin w_n)y \equiv c_n x - s_n y,\quad 
\tilde y = (\sin w_n)x + (\cos w_n)y \equiv s_n x + c_n y.
\end{equation*}
First observe that the summand theorem for trigonometric functions indicates that
\begin{align*}
x_i &= r\cos (iw_n) = r\cos \{(i-1)w_n + w_n\} = r\{\cos (i-1)w_n\} \cos w_n - r\{\sin (i-1)w_n\} \sin w_n\\ 
	&= x_{i-1}c_n - y_{i-1}s_n,\\
y_i &= r\sin (iw_n) = r\sin \{(i-1)w_n + w_n\} = r\{\sin (i-1)w_n \} \cos w_n + r\{\cos (i-1)w_n\} \sin w_n\\
	&= y_{i-1}c_n + x_{i-1}s_n.
\end{align*}
Using these identifications, we have
\begin{align}
\notag
\tilde x - x_i &= (\cos w_n)x - (\sin w_n)y - x_i
	= c_n x - s_n y - \{ x_{i-1}c_n - y_{i-1}s_n \} \\
\notag
	&= c_n (x - x_{i-1}) - s_n (y - y_{i-1}),\\
\notag
\tilde y - y_i &= (\sin w_n)x + (\cos w_n)y - y_i
\notag
	= s_n x + c_n y - \{ y_{i-1}c_n + x_{i-1}s_n \}\\
\label{identity-tri}
	&= s_n(x - x_{i-1}) + c_n (y - y_{y-1}).
\end{align}
Now
\begin{align*}
F(\tilde {\bf x}, h) &= (F_1(\tilde {\bf x}, h), F_2(\tilde {\bf x}, h)),\\
\frac{-1}{m}F_1(\tilde {\bf x}, h) &= \sum_{i=0}^{n-1} \frac{(\cos w_n)x - (\sin w_n)y -x_i}{\{ \{ (\cos w_n)x - (\sin w_n)y - x_i\}^2 + \{ (\sin w_n)x + (\cos w_n)y - y_i\}^2 + h^2\}^{(m+2)/2}},
\end{align*}
while
\begin{equation*}
\frac{-1}{m}F_2(\tilde {\bf x}, h) = \sum_{i=0}^{n-1} \frac{(\sin w_n)x + (\cos w_n)y -y_i}{\{ \{ (\cos w_n)x - (\sin w_n)y - x_i\}^2 + \{ (\sin w_n)x + (\cos w_n)y - y_i\}^2 + h^2\}^{(m+2)/2}}.
\end{equation*}
From (\ref{identity-tri}), we have
\begin{align*}
\{(\cos w_n)x - (\sin w_n)y -x_i\}^2 &= \{c_n (x - x_{i-1}) - s_n (y - y_{i-1})\}^2\\
	&= c_n^2 (x - x_{i-1})^2 -2c_n s_n(x - x_{i-1}) (y - y_{i-1}) + s_n^2 (y - y_{i-1})^2,\\
 \{ (\sin w_n)x + (\cos w_n)y - y_i\}^2 &= \{y_1(x - x_{i-1}) + x_1 (y - y_{y-1})\}^2 \\
 	&= s_n^2 (x - x_{i-1})^2 + 2c_n s_n (x - x_{i-1}) (y - y_{i-1}) + c_n^2 (y - y_{i-1})^2,\\
S &\equiv \{(\cos w_n)x - (\sin w_n)y -x_i\}^2 +  \{ (\sin w_n)x + (\cos w_n)y - y_i\}^2\\
	&= (c_n^2 + s_n^2)(x - x_{i-1})^2 + (s_n^2 + c_n^2)  (y - y_{i-1})^2\\
	&= (x - x_{i-1})^2 + (y - y_{i-1})^2
\end{align*}
and hence, writing $d_i^2 \equiv (x - x_i)^2 + (y - y_i)^2 + h^2$, we know that
\begin{align*}
\frac{-1}{m}F_1(\tilde {\bf x}, h) &= \sum_{i=0}^{n-1} \frac{(\cos w_n)x - (\sin w_n)y -x_i}{ d_i ^{(m+2)/2}} \\
	&=  \sum_{i=0}^{n-1} \frac{ c_n (x - x_{i-1}) - s_n (y - y_{i-1}) }{ d_{i-1}^{(m+2)/2} }\\
	&= c_n F_1({\bf x}, h) - s_n F_2({\bf x}, h), \\
\frac{-1}{m}F_2(\tilde {\bf x}, h) &= \sum_{i=0}^{n-1} \frac{(\sin w_n)x + (\cos w_n)y -x_i}{ d_i^{(m+2)/2}} \\
&= \sum_{i=0}^{n-1} \frac{ s_n (x - x_{i-1}) + c_n (y - y_{y-1}) }{ d_{i-1}^{(m+2)/2}} \\
	&= s_n F_1({\bf x}, h) + c_n F_2({\bf x}, h),
\end{align*}
where we have identified $d_{-1}$ with $d_{n-1}$.
As a consequence, we have $F(\xi_n\cdot {\bf x}, h) = \xi_n \cdot F({\bf x}, h)$ for all ${\bf x}\in \mathbb{R}^2$ and hence the equivariance is verified.
\end{proof}

The above arguments provide a strategy to find the equivariance of $F$ associated with a group $\Gamma$ with different point configurations, such as square-, or rectangular-like ones.
\par
\bigskip
From the equivariance of $F$, we can consider the bifurcation problem for $F$ with the height $h > 0$ as the bifurcation parameter as that with the symmetry $D_n$.
According to the general theory \cite{GSS1988}, symmetry breaking bifurcations {\em generically}\footnote{
Detailed meaning of this concept can be found in \cite{GSS1988}.
} occur following {\em isotropy subgroups}.
\begin{dfn}\rm
Let $\Gamma$ be a (compact Lie) group acting on a vector space $V$.
The {\em isotropy subgroup} of ${\bf x}\in V$ is
\begin{equation}
\Sigma_{\bf x} = \{\gamma \in \Gamma \mid \gamma \cdot {\bf x} = {\bf x} \}.
\end{equation}
Next, let $\Sigma \subset \Gamma$ be a subgroup.
The {\em fixed-point subspace} of $\Sigma$ is 
\begin{equation*}
{\rm Fix}(\Sigma) = \{{\bf x} \in V \mid \sigma \cdot {\bf x} = {\bf x} \text{ for all }\sigma \in \Sigma\}.
\end{equation*}
Finally, the action of $\Gamma$ on $V$ is said to be {\em absolutely irreducible} if the only linear mappings on $V$ that commute with either $\Gamma$ are scalar multiples of the identity.
\end{dfn}

The general theory in \cite{GSS1988} shows that the dihedral group $D_n$ includes $\mathbb{Z}_2$ as, up to isomorphisms,  the only (and maximal) isotropy subgroups, while the precise description of these subgroups is mentioned later. 
Maximality of an isotropic subgroup $\Sigma$ means that, if $\tilde \Sigma$ is an isotropic subgroup of $\Gamma$ properly containing $\Sigma$, then $\tilde \Sigma = \Gamma$.
\par
We can apply the {\em equivariant branching lemma} mentioned in \cite[Chapter XIII-3]{GSS1988}  stating as follows to characterizing symmetry-breaking structures in the bifurcation problem $F = F(\cdot, h) = 0$.

\begin{thm}[Equivariant branching lemma, \cite{GSS1988}]
\label{thm-branching}
Let $\Gamma$ be a Lie group acting absolutely irreducibly on a vector space $V$ and let $g = g({\bf x}, h)$ be a $\Gamma$-equivariant bifurcation problem satisfying
\begin{equation}
\label{absolute-g}
(D_{\bf x}g)_{{\bf 0}, h} = c(h) I\quad \text{ with }\quad c(h_0) = 0,\quad c'(h_0)\not = 0.
\end{equation}
Let $\Sigma$ be an isotropy subgroup satisfying 
\begin{equation*}
\dim {\rm Fix}(\Sigma) = 1.
\end{equation*}
Then there exists a unique smooth solution branch to $g = 0$ such that the isotropy subgroup of each solution is $\Sigma$.
\end{thm}

\begin{rem}\cite[Chapter I]{GS2002}
The (maximal) isotropic subgroups $\mathbb{Z}_2$ of $D_n$ acting on $\mathbb{R}^2$ are characterized as follows.
The simplest isotropic subgroup is 
\begin{equation*}
\Sigma = \mathbb{Z}_2(\kappa) \equiv \{1, \kappa\}, \kappa: \begin{pmatrix}
x \\ y
\end{pmatrix}\mapsto \begin{pmatrix}
x \\ -y
\end{pmatrix},
\end{equation*}
i.e., the reflection symmetry. 
Because we are interested in $D_n$-action on $\mathbb{R}^2$, then the fixed-point subspace of $\Sigma$ is 
\begin{equation*}
{\rm Fix}(\Sigma) = \{(x,0) \in \mathbb{R}^2\} \cong \mathbb{R}.
\end{equation*}
Therefore $\Sigma = \mathbb{Z}_2(\kappa)$ satisfies the requirement stated in Theorem \ref{thm-branching}.
In general, as mentioned in \cite{GS2002}, group related solutions have {\em conjugate} isotropy subgroups.
Indeed, if $g({\bf x}, h) = 0$ for the $\Gamma$-equivariant map $g = g(\cdot, h)$ for any fixed parameter $h$, then $g(\gamma\cdot {\bf x}, h) = \gamma g({\bf x}, h) = 0$ for any $\gamma \in \Gamma$.
Now suppose that $\sigma \cdot {\bf x} = {\bf x}$, namely $\sigma$ is an element of the isotropy subgroup $\Sigma_{\bf x}$.
Then we have
\begin{equation*}
(\gamma \circ \sigma \circ \gamma^{-1})\gamma \cdot {\bf x} = (\gamma \circ \sigma) \cdot {\bf x} = \gamma \cdot {\bf x}.
\end{equation*}
This implies that the following {\em conjugacy} relation between isotropic subgroups holds:
\begin{equation*}
\Sigma_{\gamma \cdot {\bf x}} = \gamma \Sigma_{\bf x} \gamma^{-1} \equiv \{ \gamma \circ \sigma \circ \gamma^{-1}\mid \sigma \in \Sigma_{\bf x}\}.
\end{equation*}
Also, if $T = \gamma \Sigma \gamma^{-1}$, then we have ${\rm Fix}(T) = \gamma ({\rm Fix}(\Sigma))$.
In particular, $\dim {\rm Fix}(T) = 1$ if $\dim {\rm Fix}(\Sigma) = 1$ and hence we have the conjugate class of $1$-dimensional isotropic subgroups.
\par
When $n$ is odd, any power of the rotation action $\xi_n^k$ on $D_n$ does not coincide with the reflection $\kappa$, indicating that the reflection symmetry group $\mathbb{Z}_2(\kappa)$ admits exactly $n$ conjugate groups $\{\mathbb{Z}_2(\kappa \circ \xi_n^k)\}_{k=0}^{n-1}$, all of which admits $1$-dimensional fixed-point subspaces:
\begin{equation}
\label{fix-odd}
{\rm Fix}(\mathbb{Z}_2(\kappa \circ \xi_n^k)) = \left\{ \left(x\cos \frac{k\pi}{n}, -x\sin \frac{k\pi}{n} \right) \mid  x \in \mathbb{R} \right\}.
\end{equation}
When $n$ is even, on the other hand, the reflection symmetry $\kappa$ coincides with the $\pi$-rotation symmetry $\xi_n^{n/2}\equiv R_{\pi/2}$.
Therefore the subgroups conjugate to $\mathbb{Z}_2(\kappa)$ are written by 
\begin{equation*}
\mathbb{Z}_2(\kappa \circ \xi_n^k)\quad \left(k=0,1,\ldots, \frac{n}{2}-1\right),\quad \mathbb{Z}_2(\kappa \circ R_{\pi/2} \circ \xi_n^k)\quad \left(k=0,1,\ldots, \frac{n}{2}-1\right), 
\end{equation*}
while fixed-point subspaces are described by
\begin{equation}
\label{fix-even}
\left\{ \left(x\cos \frac{k\pi}{n}, -x\sin \frac{k\pi}{n} \right) \mid  x \in \mathbb{R} \right\},
\end{equation}
for any $k = 0,\ldots, n-1$.
\end{rem}

Here we pay our attention to $\mathcal{X}_n = \{{\bf x}_i\}_{i=0}^{n-1} \subset \mathbb{R}^2$ with $n\geq 3$, the configuration of given points considered in Proposition \ref{prop-equivariant-LED}. 
Now suppose that one point ${\bf x}_i \in \mathcal{X}_n$ is on the $x$-axis.
We then observe from the above arguments that
\begin{itemize}
\item $F$ is $D_n$-equivariant generated by (\ref{Dn}) (Proposition \ref{prop-equivariant-LED});
\item ${\rm Fix}(D_n) = \{{\bf 0}\}$,
\item for any $n$, $\mathbb{Z}_2(\kappa \circ \xi_n^k)$, $k=0,\cdots, n-1$, are maximal isotropic subgroups of $D_n$ whose fixed-point subspaces are one-dimensional.
\end{itemize}

We shall observe in Section \ref{section-numerical-bifurcation} that there is a critical $h = h_0 > 0$ such that the Jacobian matrix $D_{\bf x}F({\bf 0}, h_0)$ is the zero matrix, while $D_{\bf x}F({\bf 0}, h)$ has double (in particular, identical) eigenvalues, which indicates that the action $D_n$ on $\mathbb{R}^2$ is absolutely irreducible and (\ref{absolute-g}) holds for $D_{\bf x}F({\bf 0}, h)$.
Theorem \ref{thm-branching} then provides the existence of the smooth curve of solutions bifurcated from the origin ${\bf 0} =(0,0)^T$ whose isotropic subgroups are determined by the above conjugate class of isotropic subgroups.
In particular, the following consequences hold:

\begin{itemize}
\item for any $n$, there is the locally unique smooth curve of solutions starting from the origin ${\bf 0} =(0,0)^T$ running on the $x$-axis;
\item there are also locally unique smooth curve of solutions starting from the origin ${\bf 0} =(0,0)^T$ running on the fixed-point subspaces (\ref{fix-odd}) or (\ref{fix-even}) for $k=1,\ldots, n-1$.
\end{itemize}
Note that, when $n$ is even, there are two bifurcated solutions to $F(\cdot, h) = 0$ on the fixed-point subspaces due to the coincidence of $\kappa$ and the $\pi$-rotation symmetry $R$. 
The most important point is that the above scenario is {\em the only possible} bifurcation from the origin at a bifurcation point for the evenly distributed circular LED configuration.
Note again that the above structure is intrinsically determined by {\em isotropic subgroups} of a Lie group $\Gamma$ acting on $\mathbb{R}^2$ and the functional $F$ being $\Gamma$-equivariant.
The similar arguments to the above ones will provide the corresponding bifurcation structures of maximizers in the other configurations.

Finally note that, as far as the bifurcation structure for $D_n$-symmetry is generated following Theorem \ref{thm-branching}, there is {\em no} possibility of the presence of bifurcations with rotational symmetry breaking.
Indeed, $D_n$ does not contain the cyclic group $\mathbb{Z}_n$ (describing rotation symmetry) as an isotropy subgroup except the origin.
Even if the isotropy subgroup including rotational symmetry, which is $\Sigma = \Sigma_{{\bf 0}} = D_n$, the fixed-point subspace is $0$-dimensional, which is incompatible with the situation in Theorem \ref{thm-branching}.

\section{Numerical Simulations}
\label{num}

\subsection{Simulations for A Realistic Work Office Model}
\label{real-office}

To illustrate the aggregated optical wirelessly transmitted power received from an LED array, as shown in Eq. \eqref{gain-total} and Eq. \eqref{gain-total-1}, simulations are carried out in the following. A realistic work office with $2.5$ m wide and $2.5$ m long is equipped with a ceiling LED array containing of $4$ rows, each row consists of $112$ LEDs, is considered for the simulation conditions. The LED array is at the center of the room. The distance between two LEDs in a row is 1 cm, between two rows in a pair is $0.01$ m, and between two row pairs is $0.1$ m. Suppose that each LED has an output power of $100$ mW, semi-angle at half illuminance of $70$ degrees. On the other hand, the optical receiver has a physical area of $2.5$ cm$^2$ and an FOV of $60$ degrees. Furthermore, the planes containing the LED array and the solar cell are in parallel. 

\begin{figure}[htpb!]
		\centering
		\includegraphics[scale=0.5]{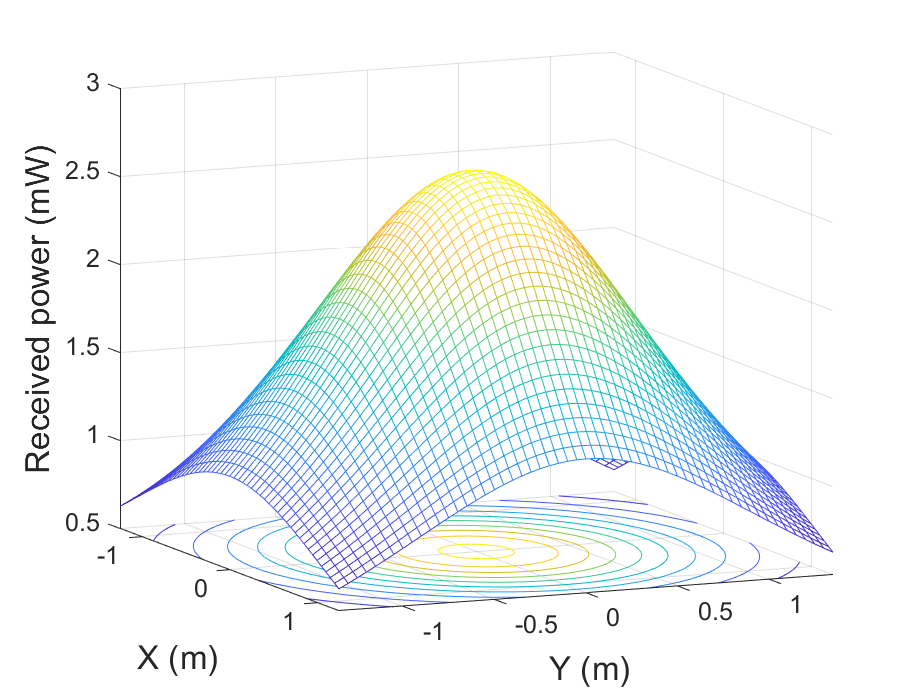}
		\caption{Wirelessly transmitted power by an LED array at a distance of $1.6$ m inside a $2.5$m $\times$ $2.5$m work office.}
		\label{LED_array_los_160}
	\end{figure}

Then the first OWPT simulation result is displayed in Fig. \ref{LED_array_los_160}, showing how the transmitted power is distributed in the horizontal plane located $1.6$ m far from the LED array. Particularly, the wirelessly transmitted power is maximum when the solar cell is perfectly aligned with the center of the LED array. As long as the solar cell is moved away from that center, the power it can receive is reduced following parabolic curves. 

\begin{figure}[htpb!]
		\centering
		\includegraphics[scale=0.5]{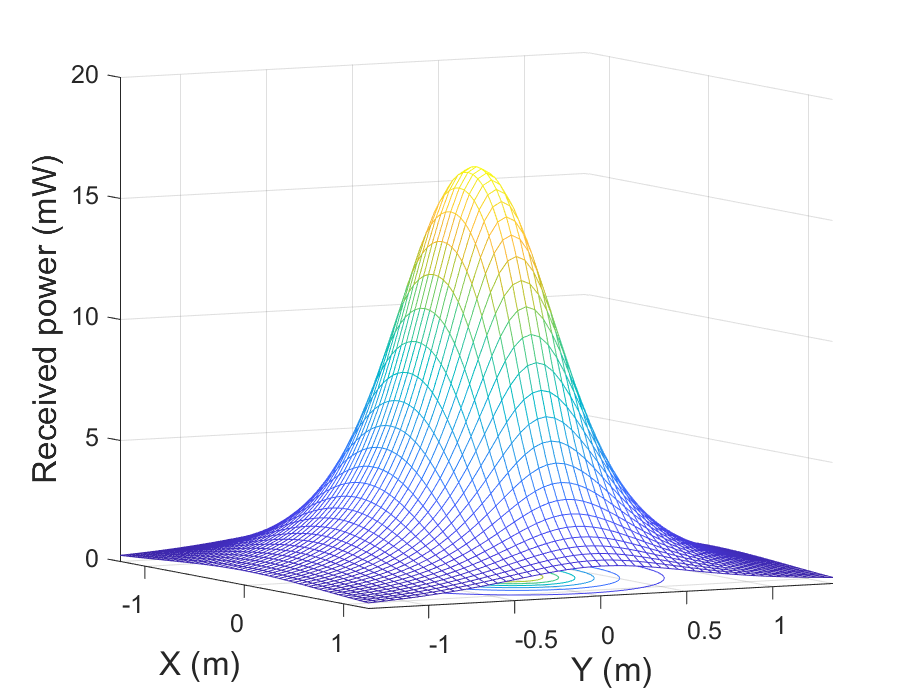}
		\caption{Wirelessly transmitted power by an LED array at a distance of $0.52$ m inside a $2.5$m $\times$ $2.5$m work office.}
		\label{LED_array_los_52}
	\end{figure}

In the second simulation, the transmitting distance is reduced to be $0.52$ m, whilst all other parameters are kept the same. The simulation result in this case is exhibited in Fig. \ref{LED_array_los_52}. Obviously, a much steeper decrease of transmitted OWPT power can be observed as the receiver moves further from the center, compared to that in Fig. \ref{LED_array_los_160}, due to the wider angle of incidence and the squared inverse dependence on the transmitting distance.

\begin{figure}[htpb!]
		\centering
		\includegraphics[scale=0.5]{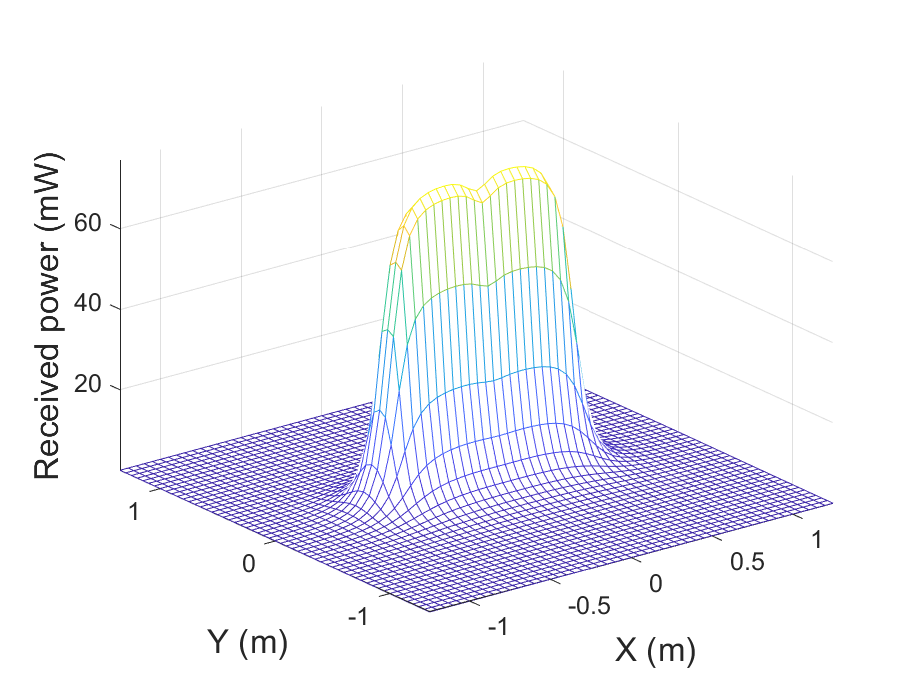}
		\caption{Wirelessly transmitted power by an LED array at a distance of 10 cm inside a 2.5m $\times$ 2.5m work office.}
		\label{LED_array_los_10}
	\end{figure}

However, when the transmitting distance is small, compared to the length of the LED array, the peak wirelessly transmitted power is not obtained at the position aligned with the LED array center, as depicted in the simulation result in Fig. \ref{LED_array_los_10}. In fact, there are two maximum power points which are are equal and symmetric through the LED array center, as can be seen in Fig. \ref{LED_array_los_10}. Additionally, OWPT is only available in a small area around the LED array. To this end, all simulations depict that OWPT for small devices is feasible with existing LED lighting systems, even at a long transmitting distance or with a large angle of incidence or irradiation.

\begin{figure}[htpb!]
		\centering
		\includegraphics[scale=0.5]{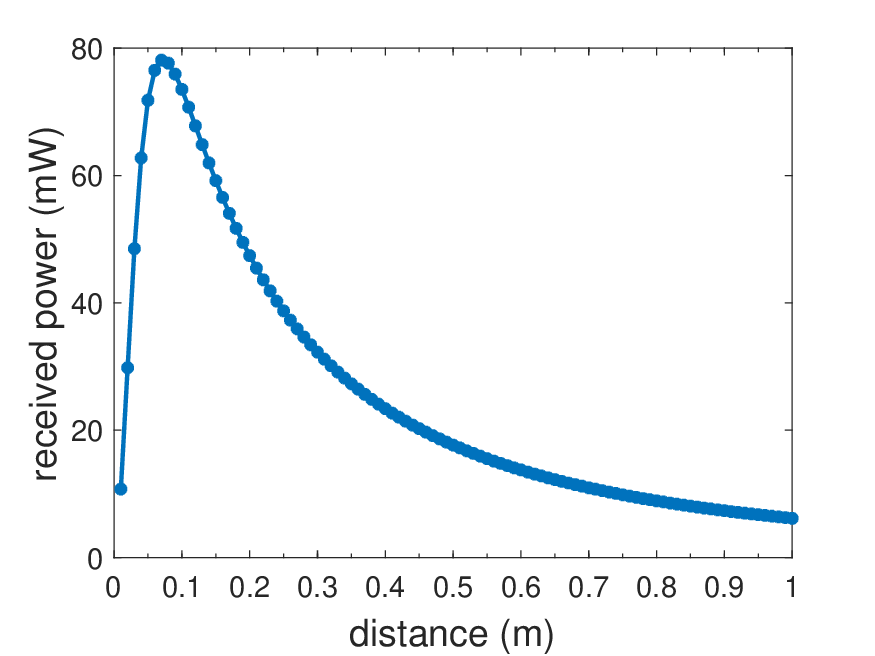}
		\caption{Variation of the wirelessly transmitted power to a receiver as its distance to the LED array is changed.}
		\label{LED_array_h_change}
	\end{figure}

Finally, we assess how the wirelessly transmitted power at a point with fixed $x$ and $y$ coordinates is varied when its $z$ coordinate, i.e., the distance between the LED array and the receiver, is changed. The simulation result displayed in Fig. \ref{LED_array_h_change} shows such variation of wirelessly transmitted power as the receiver is put perpendicular to the center of the LED array. Interestingly, the power curve exhibited in Fig. \ref{LED_array_h_change} reveals that wirelessly transmitted power is not a monotonic function of the transmitting distance. Instead, there exists a specific distance up to which the wirelessly transmitted power is monotonically increasing and beyond that the wirelessly transmitted power is monotonically decreasing. In other words, there exists a distance $h$ from the center of the LED array at which the wirelessly transmitted power is maximum over all other distances from that center. For the considering example, such distance is about $0.07$ m. 

\subsection{Simulations for Points on A Circle}
\label{unreal-num}

In this section, we illustrate the complexity of the considered optimization problem presented in Problem \ref{prob-def} for a synthetic scenario in which there are $20$ LEDs with same parameters as that in Section \ref{real-office} but now are located on a circle centered at the origin with a radius of $1.2$ m. For the purpose of comparison with the results in Section \ref{real-office}, the same receiver plane of $2.5$m $\times$ $2.5$m is considered. 

	\begin{figure}[htpb!]
		\centering
		\includegraphics[scale=0.5]{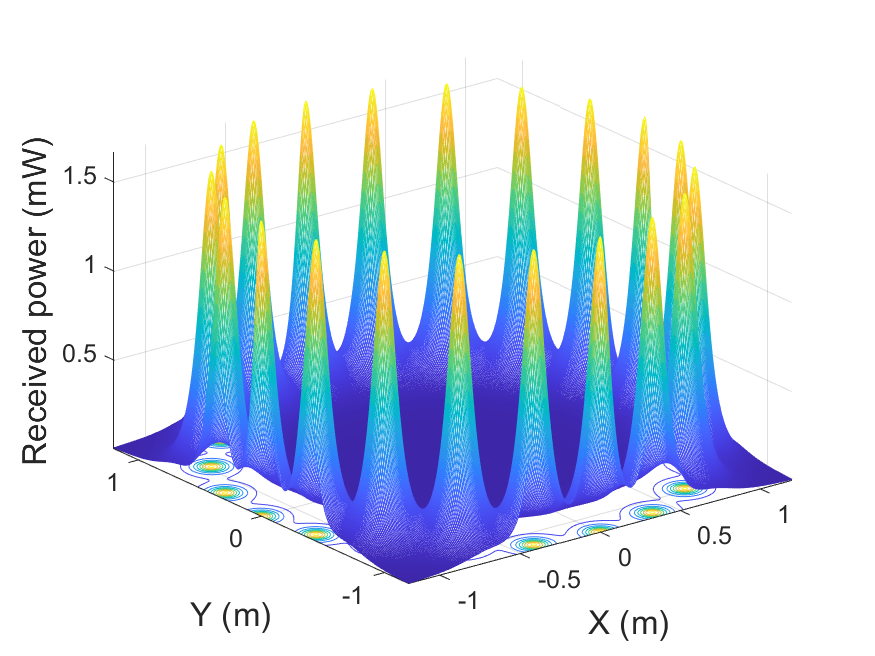}
		\caption{Distribution of the wirelessly transmitted power to a receiver distanced at $0.1$ m, with $20$ LEDs evenly located on a circle with center at the origin and radius of $1.2$ m.}
		\label{LED-r12-n20-h01}
	\end{figure}	
	
First, the simulation result for a distance of $0.1$ m between the receiver and the LEDs is shown in Fig. \ref{LED-r12-n20-h01}. It can be observed that there are $20$ global maximizers and many local optimizers in this context. Interestingly, the number of global maximizers in this case is equal to the number of LEDs. The distribution of wirelessly transmitted power is much more complex than that in Fig. \ref{LED_array_los_10} with the same transmitting distance.

\begin{figure}[htpb!]
		\centering
		\includegraphics[scale=0.5]{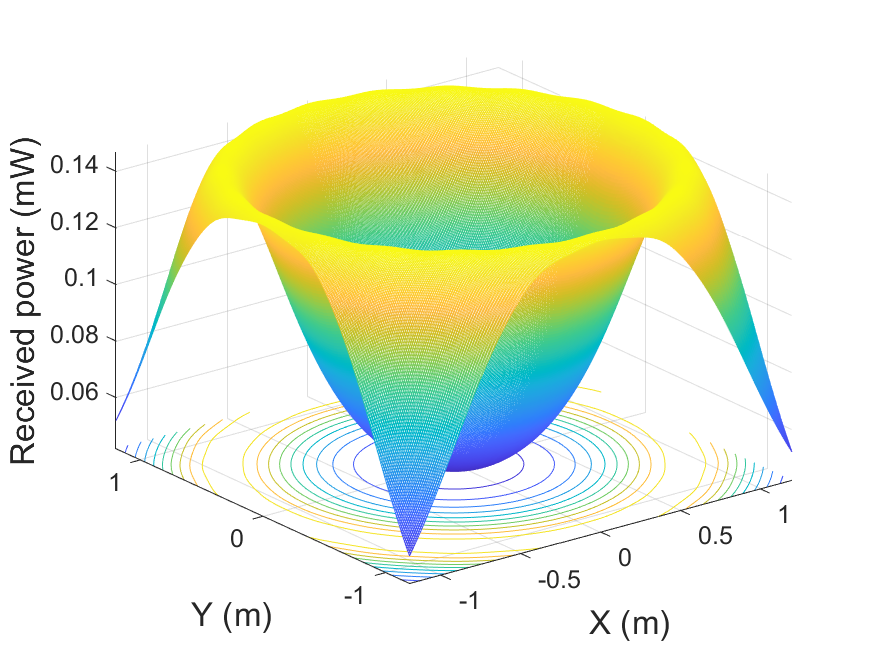}
		\caption{Distribution of the wirelessly transmitted power to a receiver distanced at $0.52$ m, with $20$ LEDs evenly located on a circle with center at the origin and radius of $1.2$ m.}
		\label{LED-r12-n20-h052}
	\end{figure}

\begin{figure}[htpb!]
		\centering
		\includegraphics[scale=0.5]{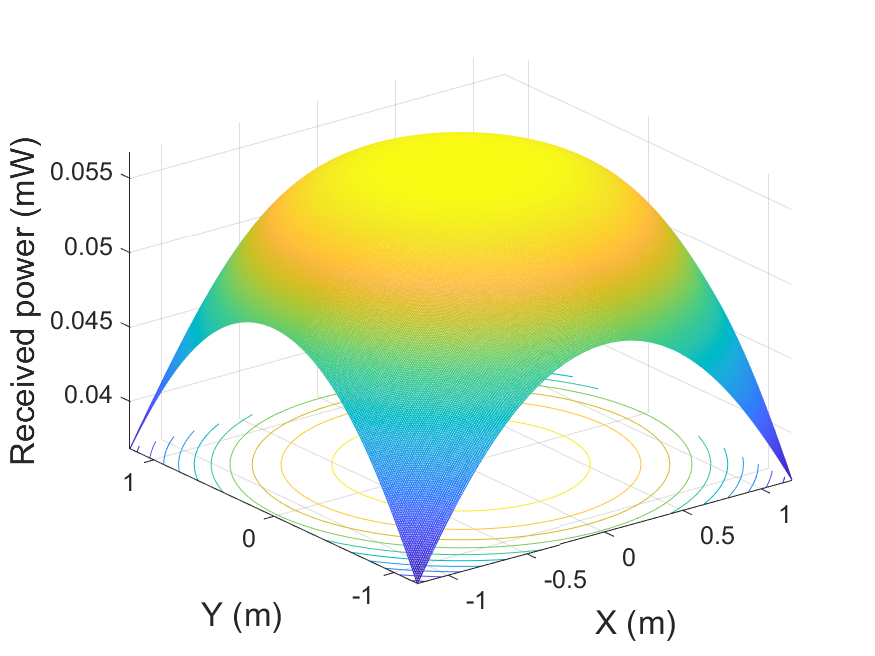}
		\caption{Distribution of the wirelessly transmitted power to a receiver distanced at $1.6$ m, with $20$ LEDs evenly located on a circle with center at the origin and radius of $1.2$ m.}
		\label{LED-r12-n20-h16}
	\end{figure}

Next, we increase the transmitting distance to be $0.52$ m and $1.6$ m, and display the simulation results in Figures \ref{LED-r12-n20-h052} and \ref{LED-r12-n20-h16}, respectively. As can be seen, the wirelessly transmitted power distribution is significantly changed as the transmitting distance is increased. Furthermore, there are many global optimizers, not a single one, in Figures \ref{LED-r12-n20-h052} and \ref{LED-r12-n20-h16}. Our simulation shows that the switch from the existence of many global optimizers to one global optimizer occurs at the distance of $1.621$ m. 
	

\subsection{Bifurcation visualization}
\label{section-numerical-bifurcation}

Here $h$-dependent families of (local) maximizers are demonstrated for several LED configurations.
Two examples here follow and extend results to detect maximizers for landscapes of the energy functional $f$ to various heights, while the rest demonstrates the tracking of maximizers for random LED configurations for showing an effectiveness of our approach in such cases.
Bifurcation diagrams as objects tracking $h$-dependent collection of are computed through {\em numerical bifurcation theory} (e.g., \cite{DKK1991_finite, K2013}).
See references for details.

\subsubsection{Lattice configurations}

\begin{figure}[htbp]\em
\centering
\includegraphics[width=.6\textwidth]{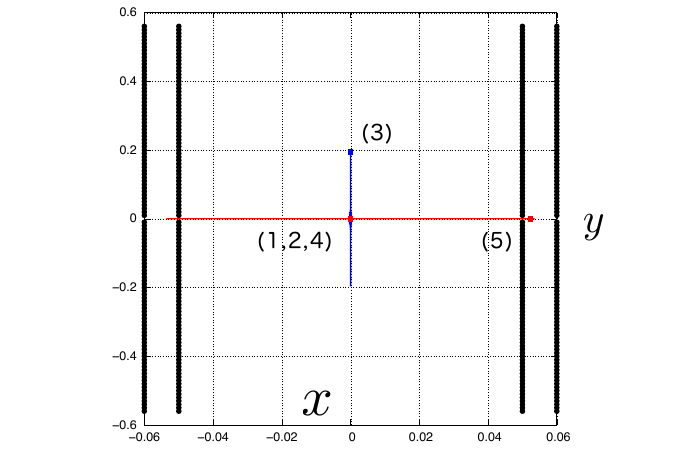}
\caption{$(x,y)$-projection of Figure \ref{fig:align_2d}}
\flushleft
\label{fig:align_2d}
Black dots represent LED locations.
Colors of branches and labels correspond to those in Figure \ref{fig:align_3d}.
\end{figure}

First the lattice configuration of LEDs consisting of
\begin{itemize}
\item $4$ arrays of LEDs whose distance between adjacent arrays is $0.1$ meters with the symmetry along the $y$-axis, while
\item each array consisting of $112$ LEDs whose distance between adjacent LEDs is $0.01$ meter except across the $x$-axis, where the adjacent row arrays have $0.02$ meters distance,
\end{itemize}
is considered (see Figure \ref{fig:align_2d}).
We can show through the similar arguments to Proposition \ref{prop-equivariant-LED} that the functional $F$ with the present configuration is equivariant under the $\Gamma = \mathbb{Z}_2\oplus \mathbb{Z}_2$-action on $\mathbb{R}^2$, where the corresponding symmetries are $\kappa_x: (x,y)^T\mapsto (-x, y)^T$ and $\kappa_y: (x,y)^T\mapsto (x, -y)^T$.
In particular, for fixed $h>0$, if ${\bf x}$ is a zero of $F(\cdot, h)$, then so are $\kappa_x({\bf x})$ and $\kappa_y({\bf x})$.
The maximal (proper) isotropy subgroups of $\Gamma$ are $\mathbb{Z}_2 = \mathbb{Z}_2(\kappa_x)$ and $\mathbb{Z}_2(\kappa_y)$, and the (generic) symmetry breaking bifurcations are followed by Theorem \ref{thm-branching}.
\par
\begin{figure}[htbp]\em
\centering
\includegraphics[width=.6\textwidth]{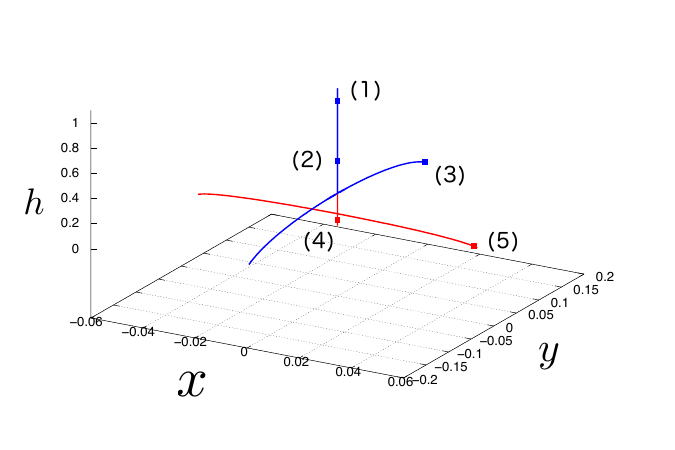}
\caption{Collection of critical points (maximizers) in a lattice LED configuration}
\flushleft
\label{fig:align_3d}
Branches of critical points as $h$-parameter families are drawn.
Blue curves correspond to critical points with negative eigenvalues, while red ones to critical points admitting at least one positive eigenvalue.
Detailed information of points with numbered labels are provided in Table \ref{table-align}.
\end{figure}

Computations of the collection of critical points drawn in Figures \ref{fig:align_3d} show that, nontrivial critical points are bifurcated from the origin ${\bf x} = (0,0)^T\in \mathbb{R}^2$ at $h \approx 0.265718574$ (meter), where the origin admits saddle-like energy landscape as observed in Figure \ref{LED_array_los_10}, and bifurcated solutions are invariant under $\mathbb{Z}_2(\kappa_x)$.
On the other hand, continuation of the origin ${\bf x} = (0,0)^T$ also yields that there is another bifurcation of critical points through exchange of signs of remaining (negative) eigenvalues and saddle-like critical points are bifurcated.
The height of bifurcation is $h \approx 0.10313732$. 
In any case, one confirms that our bifurcation diagrams indeed obey symmetries mentioned before.
\par
Table \ref{table-align} provides concrete information of critical points such as local geometric characterization of energy landscape as well as corresponding values of the power $f$, which are compatible with observations in Section \ref{real-office}.

\begin{table}[ht]\em
\centering
{
\begin{tabular}{ccccc}
\hline
Label & $h$ & ${\bf x} = (x,y)$ & $f({\bf x})$ & $\sharp$\text{positive eigenvalues}\\
\hline\\[-2mm]
(1) & $1.0$ & $(0.0, 0.0)$ & $6.16845653$ & $0$ \\ [1mm]
(2) & $0.5202552$ & $(0.0, 0.0)$ & $16.7598622$ & $0$ \\ [1mm]
(3) & $0.1104130$ & $(0.0, 0.19453106)$ & $0.39715870$ & $0$ \\ [1mm]
(4) & $0.05068238$ & $(0.0, 0.0)$ & $72.2702944$ & $2$ \\ [1mm]
(5) & $0.05171632$ & $(0.05234554, 0.0)$ & $104.470261$ & $1$ \\ [1mm]
\hline 
\end{tabular}%
}
\caption{Location of maximizers and the values of functionals in an aligned configuration}
\label{table-align}
\flushleft
\lq\lq Label" corresponds to that in Figures \ref{fig:align_2d}--\ref{fig:align_3d}.
\lq\lq $\sharp$positive eigenvalues" corresponds to local geometric features of $f$; \lq\lq $0$" corresponds to be concave, while \lq\lq $1$" to be saddle-like and \lq\lq $2$" to be convex.
\end{table}%

\subsubsection{Random configurations}

\begin{figure}[htbp]\em
\centering
\includegraphics[width=.8\textwidth]{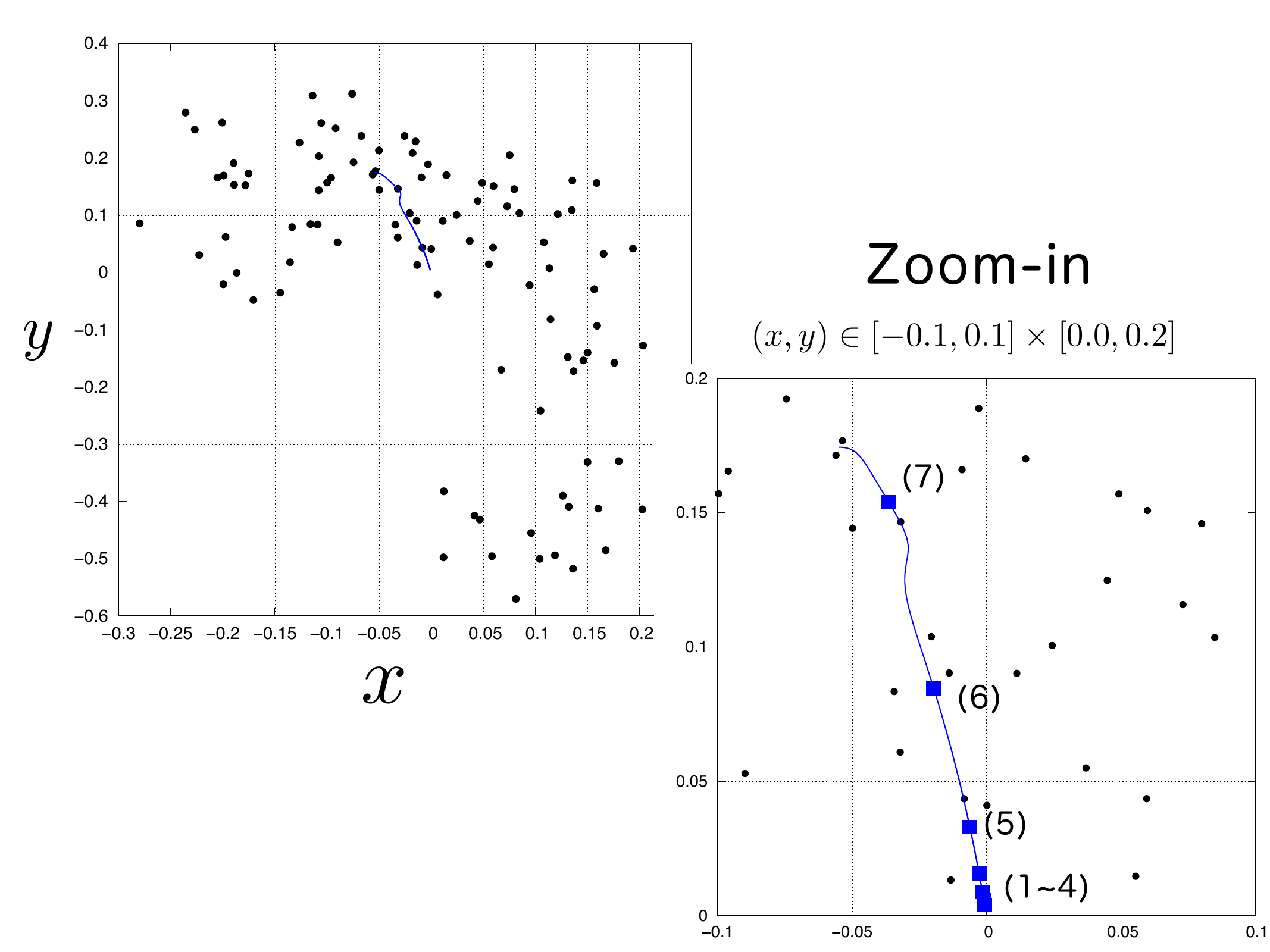}
\caption{$(x,y)$-projection of Figure \ref{fig:random_2d}}
\flushleft
\label{fig:random_2d}
Black dots represent LED locations.
Colors of branches and labels correspond to those in Figure \ref{fig:random_3d}.
\end{figure}

\begin{figure}[htbp]\em
\centering
\includegraphics[width=.6\textwidth]{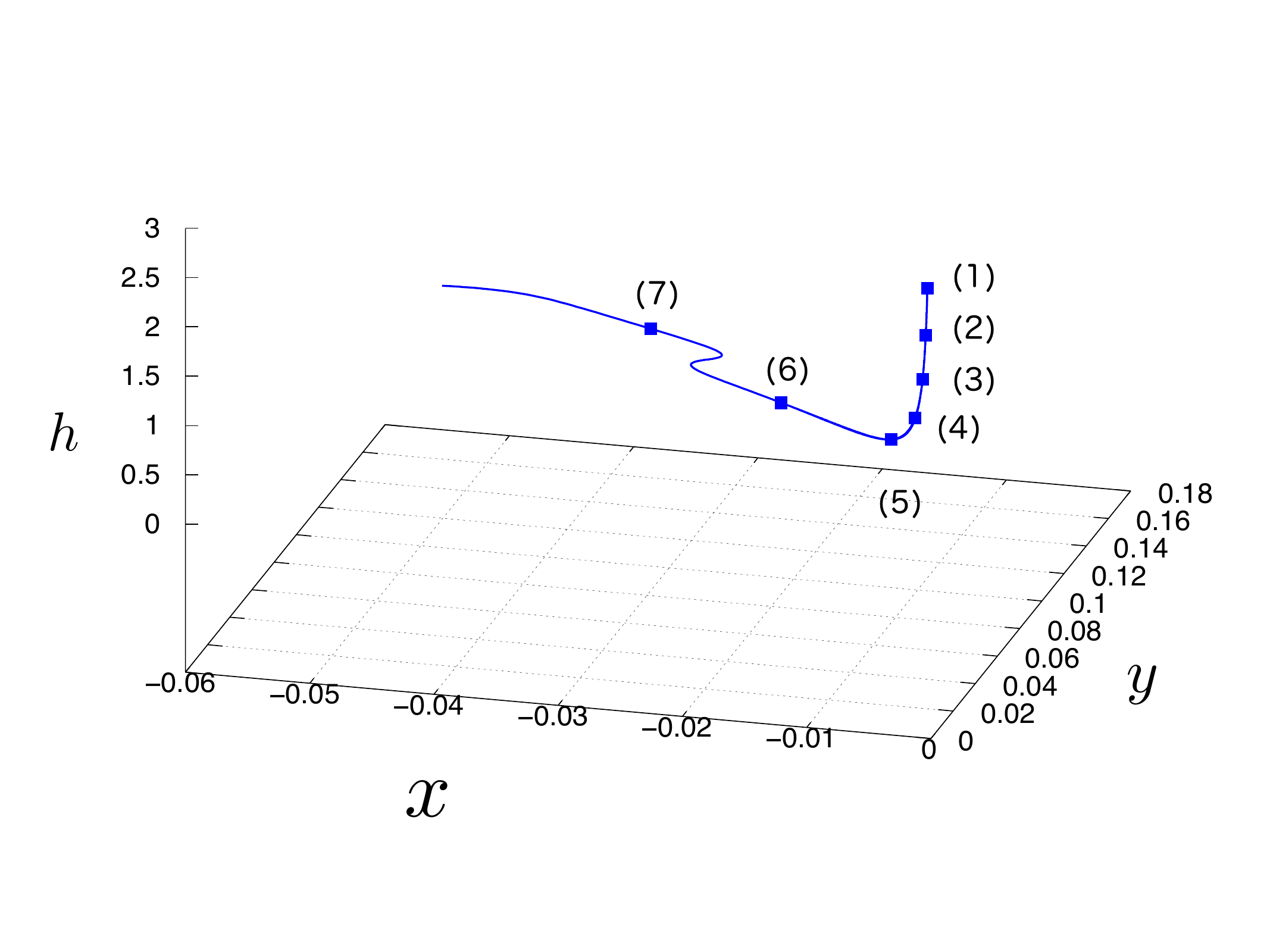}
\caption{Collection of critical points (maximizers) in a randomly distributed LED configuration}
\flushleft
\label{fig:random_3d}
Branches of critical points as $h$-parameter families are drawn.
Blue curve corresponds to the collection of critical points with negative eigenvalues, which does not admit any bifurcations as far as we have calculated.
\end{figure}

Next, consider a randomly distributed LED configuration.
Figures \ref{fig:random_2d}--\ref{fig:random_3d} show the sample configuration with $100$ LEDs and the collection of maximizers of $f$ depending on $h$.
Starting our computations at $h = 3.0$, bifurcation algorithm \cite{DKK1991_finite} works well to extract critical points with strictly negative definite Jacobian matrices at points on the curve, indicating that the collection of critical points represents maximizers of $f$ indeed.
Details on the critical points and the functional values at those points are given in Table \ref{table-random}.

\begin{table}[ht]\em
\centering
{
\begin{tabular}{cccc}
\hline
Label & $h$ & ${\bf x} = (x,y)$ & $f({\bf x})$ \\
\hline\\[-2mm]
(1) & $3.0$ & $(-0.0006372, 0.0039683)$ & $0.17945693$\\ [1mm]
(2) & $2.5$ & $(-0.0009241, 0.0056995)$ & $0.25689953$\\ [1mm]
(3) & $2.0$ & $(-0.0146160, 0.0088581)$ & $0.39715870$\\ [1mm]
(4) & $1.5$ & $(-0.0026545, 0.0155297)$ & $0.69068483$\\ [1mm]
(5) & $1.0$ & $(-0.0061445, 0.0330292)$ & $1.46853717$\\ [1mm]
(6) & $0.5$ & $(-0.0196216, 0.0847345)$ & $4.77352797$\\ [1mm]
(7) & $0.1$ & $(-0.0363161, 0.1539201)$ & $30.2078349$\\ [1mm]
\hline 
\end{tabular}%
}
\caption{Location of maximizers and the values of functionals in a random configuration}
\label{table-random}
\flushleft
\lq\lq Label" corresponds to that in Figures \ref{fig:random_2d}--\ref{fig:random_3d}.
\end{table}%

\subsubsection{Circular configurations: symmetry and possible maximizers of $f$ with relatively large heights}

Next the circular configuration with $20$ LEDs evenly located on the circle with center at the origin and radius of $1.2$ m is considered (cf. Figures \ref{LED-r12-n20-h052} and \ref{LED-r12-n20-h16}).
As observed in Section \ref{section-symmetry}, the corresponding functional $F$ is $D_{20}$-equivariant, where $D_{20}$ is the dihedral group generated by rotations and reflections given in (\ref{Dn}) with $n=20$.
Tracking critical points, we see that eigenstructure of the Jacobian matrix origin changes at $h = h_0 \approx 1.622$.
More precisely, the Jacobian matrix at ${\bf x} \in \mathbb{R}^2$ with fixed $h$ is given as follows:
\begin{align*}
D_{\bf x}F({\bf x}, h) &= \begin{pmatrix}
a_{11} & a_{12} \\ a_{21} & a_{22}
\end{pmatrix},\\
a_{11} &= \frac{\partial F_1}{\partial x}({\bf x}, h) = -m \sum_{i=1}^n \frac{1}{d_i^{m+2}} \left\{ 1 - (m+2) \frac{(x-x_i)^2}{d_i^{m+4}}\right\},\\
a_{12} &= \frac{\partial F_1}{\partial y}({\bf x}, h) = \frac{\partial F_2}{\partial x}({\bf x}, h) (= a_{21}) \\
	&= m(m+2) \sum_{i=1}^n \frac{(x-x_i)(y-y_i)}{d_i^{m+4}},\\
a_{22} &= \frac{\partial F_2}{\partial y}({\bf x}, h) = -m \sum_{i=1}^n \frac{1}{d_i^{m+2}} \left\{ 1 - (m+2) \frac{(y-y_i)^2}{d_i^{m+4}}\right\}.
\end{align*}
Substituting ${\bf x} = {\bf 0}$, let
\begin{align*}
D_{\bf x}F({\bf 0}, h) = \begin{pmatrix}
a_{11}^\ast(h) & a_{12}^\ast(h) \\ a_{21}^\ast(h) & a_{22}^\ast(h)
\end{pmatrix}
\end{align*}
We then have $a_{12}^\ast(h) \equiv 0$ for all $h > 0$ by symmetry of LED configurations and $n$ being even, while $a_{11}^\ast(h) = a_{22}^\ast(h)$ for all $h>0$.
At $h = h_0$, we further have $a_{11}^\ast(h_0) = a_{22}^\ast(h_0) = 0$.
On the other hand, 
\begin{align*}
\frac{\partial^2 F_1}{\partial x\partial h}({\bf x}, h) &= m(m+2) h \sum_{i=1}^n \left\{ \frac{1}{d_i^{m+4}} - (m+4)\frac{(x-x_i)^2}{d_i^{m+6}} \right\},\\
\frac{\partial^2 F_2}{\partial y\partial h}({\bf x}, h) &= m(m+2) h \sum_{i=1}^n \left\{ \frac{1}{d_i^{m+4}} - (m+4)\frac{(y-y_i)^2}{d_i^{m+6}} \right\}
\end{align*}
and numerical computations yield
\begin{equation*}
\frac{\partial^2 F_1}{\partial x\partial h}({\bf 0}, h_0) = \frac{\partial^2 F_2}{\partial y\partial h}({\bf 0}, h_0) \approx -1.09835 \not = 0.
\end{equation*}
These observations indicate that absolute irreducibility mentioned in Theorem \ref{thm-branching} is satisfied with $V = \mathbb{R}^2$ and 
\begin{equation*}
c(h) = \frac{\partial F_1}{\partial x}({\bf 0}, h) = \frac{\partial F_2}{\partial y}({\bf 0}, h).
\end{equation*}
In particular, $n$ (maximal) isotropy subgroups $\{\mathbb{Z}_2(\kappa_y\circ \xi_n^k)\}_{k=0}^{19}$ of $D_{20}$ generate nontrivial critical points, which correspond to $20$ mountain-like energy landscapes observed in Figure \ref{LED-r12-n20-h052}.

\begin{rem}
We omit the concrete bifurcation diagram in the present example because the bifurcation is highly degenerate in the sense of the standard machinery in the numerical bifurcation theory \cite{DKK1991_finite}.
\end{rem}

\section{Conclusion}
\label{concl}

In this study, a strongly nonlinear and non-convex geometric optimization problem is investigated, which is motivated from a realistic optical wireless power transfer system. In this problem, a sum of negatively and fractionally powered distances between a freely moving point in a plane and a set of given points arbitrarily located in another plane is maximized. The aim is then to find the positions of the freely moving point in its plane such that the sum described above is maximum, i.e., to characterize the optimizer set. Accordingly, it is shown that there is a critical value that a unique global optimizer exists when the inter-plane distance is greater than it. On the other hand, when the inter-plane distance is not larger than such critical value, the bifurcation theory is employed to derive the exact number of local maximizers which is equal to the number of bifurcation branches determined via one-dimensional isotropic subgroups of a Lie group acting on $\mathbb{R}^2$. 
Extensive numerical simulations and bifurcation points calculation are then conducted under different configurations of the given points, revealing the correctness and effectiveness of the derived theoretical results.

\section*{Acknowledgements}
Dinh Hoa Nguyen was supported in part by JSPS Kakenhi Grant Number JP23K03906. 
Kaname Matsue was partially supported by World Premier International Research Center Initiative (WPI), Ministry of Education, Culture, Sports, Science and Technology (MEXT), Japan, and JSPS Grant-in-Aid for Scientists (B) (No. JP23K20813).


\bibliographystyle{plain}
\bibliography{References}

\end{document}